\documentclass[12pt]{article}

\usepackage{amsmath,amscd,amsthm,makeidx,graphicx,url,bm,bbm,mathrsfs,palatino,bm,mathtools,latexsym,hyperref}
\usepackage[bbgreekl]{mathbbol}

\usepackage{amssymb}
\usepackage{a4wide,xy,graphicx}
\usepackage{tikz-cd}
\usepackage{enumerate,comment}
\usepackage{titlesec}
\usepackage{tikz-cd}

\numberwithin{equation}{section}

\newtheorem{theorem}{Theorem}

\newtheorem{conjecture}[theorem]{Conjecture}

\theoremstyle{remark}
\newtheorem{remark}[theorem]{Remark}

\newtheorem{definition}[theorem]{Definition}

\numberwithin{theorem}{section}

\newcounter{margin}
	{\end{itshape}  \bigskip}

\def\beq{\begin{eqnarray}}
\def\eeq{\end{eqnarray}}
\def\bes{\begin{eqnarray*}}
	\def\ees{\end{eqnarray*}}

\def\calE{{\mathcal E}}

\DeclareMathOperator{\Spec}{Spec} 
\DeclareMathOperator{\Hom}{Hom}

\DeclareMathOperator{\Ext}{Ext }

\def\bE{{\mathbb{E}}}
\def\bPhi{{\mathbb{\Phi}}}

\def\A{{\mathcal A}}

\def\E{\mathbb{E}}

\def\C{\mathbb{C}}

\def\M{{\mathcal{M}}}

\def\calA{{\mathcal{A}}}

\def\calF{{\mathcal{F}}}

\def\calH{\mathcal{H}}

\def\P{\mathbb{P}}
\def\N{\mathbb{Z}_{\geq 0}}

\newcommand{\T}{{\mathbb{T}}}

\def\calC{{\mathcal C}}
\def\calG{{\mathcal G}}

\def\calO{{\mathcal O}}

\def\Z{\mathbb{Z}}

\def\gl{{\mathfrak g\mathfrak l}}

\newcommand{\nc}{\newcommand}

\newcommand{\g}{\mathfrak{g}}

\usepackage{leftidx}

\newcommand{\ad}{\textnormal{ad}}
\newcommand{\Gr}{\textnormal{Gr}}
\newcommand{\oGr}{\overline{\textnormal{Gr}}}

\newcommand{\End}{\textnormal{End}}

\nc{\op}[1]{\mathop{\mathchoice{\mbox{\rm #1}}{\mbox{\rm #1}}
		{\mbox{\rm riptsize #1}}{\mbox{\rm \tiny #1}}}\nolimits}
\nc{\al}{\alpha}

\nc{\ep}{\varepsilon} 
\nc{\ga}{\gamma} 
\nc{\Ga}{\Gamma}
\nc{\la}{\lambda} 
\nc{\La}{\Lambda} 
\nc{\si}{\sigma}
\nc{\Sig}{{\Gamma}} 
\nc{\Om}{\Omega} 
\nc{\om}{\omega}

\nc{\SL}{\mathrm{SL}} 
\nc{\GL}{\mathrm{GL}} 
\nc{\SO}{\mathrm{SO}} 
\nc{\PGL}{\mathrm{PGL}}
\newcommand{\G}{\mathrm{G}}

\nc{\W}{\mathrm{W}}
\nc{\Lg}{\mathrm{L}}
\nc{\Pg}{\mathrm{P}}
\nc{\calL}{{\mathcal L}}
\nc{\Sym}{{\mathrm Sym}}
\nc{\h}{{\rm h}}

\nc{\Frob}{\mathrm{Frob}}
\nc{\spec}{{\mathrm Spec}}

\def\U{{\mathrm{U}}}

\nc{\rN}{\mathrm{N}}

\usepackage{longtable}

\newcommand{\Dol}{{\mathrm{Dol}}}
\nc{\cpt}{{\op{cpt}}} 
\nc{\DR}{{\mathrm{DR}}}
\nc{\B}{{\mathrm{B}}} \nc{\Triv}{\op{Triv}} \nc{\Hod}{{\op{Hod}}}
\nc{\Log}{{\op{Log}}} \nc{\Exp}{{\op{Exp}}} \nc{\Est}{E_{\op{st}}}
\nc{\Hst}{H_{\op{st}}} \nc{\Left}[1]{\hbox{$\left#1\vbox to
		10.5pt{}\right.\nulldelimiterspace=0pt \mathsurround=0pt$}}
\nc{\Right}[1]{\hbox{$\left.\vbox to
		10.5pt{}\right#1\nulldelimiterspace=0pt \mathsurround=0pt$}}
\nc{\LEFT}[1]{\hbox{$\left#1\vbox to
		15.5pt{}\right.\nulldelimiterspace=0pt \mathsurround=0pt$}}
\nc{\RIGHT}[1]{\hbox{$\left.\vbox to
		15.5pt{}\right#1\nulldelimiterspace=0pt \mathsurround=0pt$}}

\nc{\bee}{{\bf E}} 

\numberwithin{equation}{section}

\begin{document}





\title{Enhanced mirror symmetry for Langlands dual Hitchin systems}

\author{ Tam\'as Hausel
	\\ {\it IST Austria} 
	\\{\tt tamas.hausel@ist.ac.at}} 
	
\maketitle


\begin{abstract}
The first part of this paper is a survey of mathematical results on mirror symmetry phenomena between Hitchin systems for Langlands dual groups. The second part introduces and discusses multiplicity algebras of the Hitchin system on Lagrangians, and considers corresponding conjectural structures on their mirror.  
\end{abstract}


 \section{Introduction} Considering the $2$-dimensional reduction of the Yang-Mills equations in $4$ dimensions, Hitchin  \cite{hitchin} in 1987 introduced and studied the moduli space of solutions to certain self-duality equations on a Riemann surface. The moduli space turns out to have an "extremely rich geometric structure". 

The moduli space $\M$ of solutions for a complex reductive structure group $\G$ (Hitchin first considered $\G=\SL_2$),  carries a canonical hyperk\"ahler metric $g$ with complex structures $I,J$ and $K$, and corresponding K\"aher forms $\omega_I$,$\omega_J$ and $\omega_K$. In complex structure $I$ it agrees with  the moduli space $\M_\Dol$ of Higgs bundles - or Hitchin pairs - $(E,\Phi)$, where $E$ is a $\G$-bundle and  the {Higgs field} $\Phi\in H^0(C;\ad(E) \otimes K)$ is a section of the adjoint bundle twisted by the canonical bundle on a complex curve $C$: $$(\M,I)\cong \M_{\Dol}.$$  Under an isomorphism induced by multiplying the Higgs field with $i$ the K\"ahler manifolds $(\M,J,\omega_J)$ and $(\M,K,\omega_K)$ are isomorphic. In turn, they are  both isomorphic with the moduli space $\M_\DR$ of flat  $\G$-connections on the curve $C$: \beq \label{JK}(\M,J)\cong (\M,K)\cong \M_{\DR}.\eeq The notation $\M_\Dol$ for Dolbeault and $\M_\DR$ for De Rham non-abelian cohomologies follows \cite{simpson} who introduced the viewpoint of non-abelian Hodge theory in the study of $\M$. 

In turn, Hitchin \cite{hitchin-stable} introduced the Hitchin map $$h:\M\to \A.$$ In the $\G=\GL_n$ case this is just the  characteristic polynomial of the Higgs field $$h(E,\Phi)=\det(x-\Phi)\in \calA:=\bigoplus_{i=1}^n H^0(C;K^i).$$ For general $\G$ one needs to consider invariant polynomials on $\g:=Lie(\G)$ and compute them on the Higgs field. He proved in \cite{hitchin-stable} that $h$ is an algebraically completely integrable Hamiltonian system with respect to the $I$-holomorphic symplectic form $\omega_\C:=\omega_J+i\omega_K$.
Thus $h$ is sometimes referred to as the {\em Hitchin system}. This means that $\dim(\A)=\dim(\M)/2$ and that the component functions of $h$ are independent and Poisson commute.  Additionally, the Hitchin map is proper, which was proved by Hitchin for $\SL_2$ in \cite{hitchin-stable}, for $\GL_n$ by Nitsure \cite{nitsure} and Simpson \cite{simpson} and by Faltings \cite{faltings} for general $\G$. The complete integrability and properness of the Hitchin map together imply that its generic fiber is a torsor over an Abelian variety. In particular topologically they are isomorphic to compact tori. 

Due to the flexibility of their constructions (choice of curves and structure groups - but also various types of ramification data) Hitchin systems have been related to most of the known integrable systems \cite{donagi-markman,donagi}. They thus play a central role in the field of integrable systems. 

Our main interest in this survey will be  how mirror symmetry and Langlands duality relate to the Hitchin system. 
In 2003 the paper \cite{hausel-thaddeus} mathematically related the Hitchin system for Langlands dual groups to mirror symmetry. In particular it formulated a topological mirror symmetry conjecture for certain $\SL_n$ and $\PGL_n$ Hitchin systems. 

In 2007 Kapustin-Witten \cite{kapustin-witten} placed the Hitchin system in the framework of a certain supersymmetric 4-dimensional Yang-Mills theory reduced to  2 dimensions. It also offered a detailed understanding of mirror symmetry and the geometrical Langlands program as a reduction of $S$-duality in $4$ dimensions. 
This lead to many papers, such as \cite{baraglia-schaposnik, garcia-prada-etal, hitchin-characteristic,hausel-mellit-pei, hausel-hitchin} discussing pairs of mirror branes in Langlands dual Hitchin systems. The last two papers emphasised a further structure on $\M$, namely a canonical $\T$-action given by $(E,\Phi)\mapsto (E,\lambda\Phi)$. 

In 2010 Ng\^{o} \cite{ngo} proved  the Fundamental Lemma in the Langlands program, via a detailed understanding of the cohomology of certain singular fibers of the Hitchin map. In \cite{hausel-gths} similarities between the topological mirror symmetry conjecture of \cite{hausel-thaddeus} and such cohomological results of Ng\^{o} were discussed and were conjecturally related. 
In \cite{groechenig-etal1} Gr\"ochenig--Wyss--Ziegler proved the topological mirror symmetry conjecture of \cite{hausel-thaddeus} using an arithmetic $p$-adic integration technique. In turn in \cite{groechenig-etal2} the same authors managed to reprove  Ng\^{o}'s  cohomological results with these new $p$-adic techniques. More recently Maulik--Shen \cite{maulik-shen} managed to complete some of the suggestions of \cite{hausel-gths} and derived a proof of the topological mirror symmetry conjecture from Ng\^{o}'s results. 

First we will discuss some of the background to these developments, and then in \S \ref{enhanced} we will explain some unpublished results about enhanced mirror symmetry for Langlands dual groups at the tip of the nilpotent cone. 

\section{Background}
\subsection{Mirror symmetry}

Three aspects of mirror symmetry will be relevant for us: topological and homological mirror symmetry and  Strominger--Yau--Zaslow mirror symmetry.

Mirror symmetry in a nutshell relates the complex geometry of a Calabi-Yau $X$ with complex structure $I_X$ and K\"ahler, in particular symplectic,  $2$-form $\omega_X$ to the symplectic geometry of a same dimensional mirror Calabi-Yau $(Y,I_Y,\omega_Y)$. Originally \cite{greene-plesser} $3$-dimensional examples of such a correspondence appeared in string theory and mirror symmetry  as the statement that the physics of a  certain $2$-dimensional  type $A$ non-linear sigma model with target $(X,\omega_X)$ matches that of a $2$-dimensional type $B$ non-linear sigma model with target $(Y,I_Y)$. 

The study of mathematical aspects of this mirror symmetry has been one of the central subjects in modern symplectic/ complex algebraic geometry. The first mathematical aspect of the mirror relationship is the agreement of Hodge numbers $h^{n-p,q}(X)=h^{p,q}(Y)$, which we call {\em topological mirror symmetry}. 

The formulation of Kontsevich \cite{kontsevich} in 1994 of {\em homological mirror symmetry}  $$Fuk(X,\omega_X)\cong D^b(Y,I_Y)$$  - the agreement of the Fukaya category of $X$ and the derived category of coherent sheaves on $Y$ - gave a profound mathematical conjecture for what mirror symmetry should mean. 

The early 90's saw several  constructions of conjectured mirror pairs in \cite{candelas-etal1,batyrev1,batyrev-borisov}. In 1996 \cite{strominger-etal} Strominger--Yau--Zaslow  suggested a way  to construct the mirror of a Calabi-Yau $3$-fold $X$ out of the geometry of $X$. They argued that there should be fibrations of $X$ and $Y$ over the same base $B\cong S^3$: \beq \label{syz}\begin{array}{ccccc} X&&&&Y \\ &\searrow&&\swarrow& \\ &&B&& \end{array}\eeq so that the generic fibers are dual special Lagrangian $3$-tori. Here $L\subset X$  special Lagrangian  means that $L$ is Lagrangian $\omega_X|_L=0$ and additionally $Im(\Omega_X)|_L=0$ the imaginary part of the Calabi-Yau volume form vanishes on $L$. In turn then $Y$ should be thought of as the moduli space of certain objects in $Fuk(X,\omega_X)$, generically defined by special Lagrangian $3$-tori  equipped with a $U(1)$-local system.  A mathematical formulation of Strominger--Yau--Zaslow was pursued by Gross--Siebert \cite{gross-siebert1}, with many accomplishments and recent breakthroughs \cite{gross-siebert2,gross-etal}. 

Between mirror Calabi--Yau $3$-folds a complete construction of dual special Lagrangian fibrations \eqref{syz} is still missing. The Higgs bundle moduli spaces for Langlands dual groups, where the Hitchin systems will automatically give us such dual special Lagrangian fibrations, is a natural example, albeit in a geometrically different scenario from the original \cite{strominger-etal}. 

\subsection{Geometric Langlands correspondence}
In the works of \cite{drinfeld1},\cite{laumon1},\cite{beilinson-drinfeld} a geometric version of the Langlands correspondence has been proposed. Recall that the Langlands program in number theory (see \cite{gelbart,frenkel} for some introductory ideas)  for a reductive group $\G$ over a number field relates automorphic data for $\G$ - like modular forms for $\SL_2$ - with spectral data - like a Galois representation on the cohomology of an elliptic curve for $\SL_2$ - for a Langlands dual $\G^L$ group.  The conjectures can be formulated over the other kind of global field as well: the function field of a curve over a finite field. The conjectures become more tractable in this case as the algebraic geometry of curves can be efficiently used. Here we will consider the even more geometric version of this program for function fields of a curve over the complex numbers. 

Over $\C$ the Langlands dual $\G^\vee:=\G^L$ of a complex reductive group $\G$ is simple to construct. The classification of complex reductive groups is via their root datum $(X,\Phi,X^\vee,\Phi^\vee)$ - consisting of
a rank $n$ lattice $X$, a root system and coroot system $\Phi\subset X$ and $\Phi^\vee\subset X^\vee$ in the dual lattice $X^\vee$, satisfying certain properties- attached to all rank $n$  complex reductive groups $\G$. The Langlands dual of $\G$ then is the reductive group $\G^\vee$ whose root datum is the dual root datum $(X^\vee,\Phi^\vee,X,\Phi).$ For example $\GL_n^\vee\cong \GL_n$ and $\SL_n^\vee\cong \PGL_n$. 

The geometric Langlands correspondence of Beilinson-Drinfeld  \cite{beilinson-drinfeld} for a smooth projective curve $C$  proposes to construct from a $\G$-local system on $C$ - a geometric analogue of a Galois representation- a holonomic $D$-module on the moduli stack of bundles  $Bun_{\G^\vee}$ - a geometric analogue of an automorphic form. The main property of this construction is that the holonomic $D$-module must be an eigensheaf of certain Hecke operators. \cite{beilinson-drinfeld} succeed in this construction for a certain set of $\G$-local systems on $C$, the so-called opers. 

We will see below how Beilinson-Drinfeld's  picture can be understood also as enhanced mirror symmetry between the hyperk\"ahler, and thus in particular Calabi-Yau, moduli space of flat $\G$ connections on $\M_\DR$ and the moduli space of flat $\G^\vee$ connections $\M^\vee_{\DR}$. 

\subsection{SYZ mirror symmetry for Langlands dual Hitchin systems} 
The starting point of \cite{hausel-thaddeus} was the observation that for $\G=\SL_n$ and $\G^\vee=\PGL_n$ the two Hitchin systems $h:\M_\Dol\to \A$ for $\G$-Higgs bundles and $\M^\vee_\Dol\to \calA^\vee$ for $\G^\vee$ Higgs bundles have the same base $\calA\cong \calA^\vee$ and in the diagram  \beq \label{syzdol}\begin{array}{rcccl} \M_\Dol &&&&\!\!\!\!\!\!\M^\vee_\Dol \\ &\!\!\!\! {}_h \!\!\searrow&&\swarrow_{h^\vee} &  \\ &&\calA&& \end{array}\eeq the generic fibers are 
dual Abelian varieties. The fibers are holomorphic Lagrangian with respect to the holomorphic symplectic forms $\omega_\C:=\omega_J+i\omega_K$, i.e. $\omega_\C=0$, thus both $\omega_J=\omega_K=0$, along the fibers because of the complete integrability of the Hitchin systems. 

The same maps in complex structure $J$    then yield \beq \label{syzdr}\begin{array}{rcccl} \M_\DR &&&&\!\!\!\!\!\!\M^\vee_\DR \\ &\!\!\!\! {}_h \!\!\searrow&&\swarrow_{h^\vee} &  \\ &&\calA&& \end{array}\eeq where the generic fibers are special Lagrangian fibrations. This means that $\omega_J$  vanishes on the fibers, thus they are  Lagrangian.  Additionally $\omega_K$ as well as the imaginary part of the $J$-holomorphic Calabi-Yau form $(\omega_K+i \omega_I)^{2d}$ also vanish along the fibers. 

To find a mathematically testable form of mirror symmetry the paper \cite{hausel-thaddeus} formulated {\em topological mirror symmetry} between the mirror Calabi-Yau's: $\M_{\DR}$ and $\M^\vee_{\DR}$.
\subsection{Topological mirror symmetry for Langlands dual Hitchin systems}

To formulate this version of mirror symmetry from \cite{hausel-thaddeus} we will be a bit more precise about our moduli spaces. For $\SL_n$ we consider $\M:=\M_\Dol$ the moduli space of stable Hitchin pairs of rank $n$ fixed bundles $E$ on $C$ of fixed determinant line bundle of  degree $1$  with trace-free Higgs fields $\Phi\in H^0(C;\End_0(E)\otimes K)$. For $\PGL_n$ we consider the action of $\Gamma:=Jac_C[n]$, the group of order $n$ line bundles on $C$, on the $\SL_n$ moduli space $\M$ and define $\M^\vee:=\M/\Gamma$. Then $\M$ will be a smooth quasi-projective variety, while $\M^\vee$ is a quasi-projective orbifold.  For our considerations we will need extra twisting structure on the moduli spaces in the form of a gerbe $$\alpha\in H^2(\M^\vee,\U(1))\cong H^2(\M,\U(1))^\Gamma$$ which can be constructed using   the universal bundle on $\M$. 

One can then define  certain mixed Hodge numbers $$h^{p,q}(\M)=h^{p,q}(H^{p+q}_c(\M))$$ for the smooth $\M$ and $\alpha$-twisted stringy Hodge numbers for the orbifold $\M^\vee=\M/\Gamma$:  $$h_{st,\alpha}^{p,q}(\M^\vee):=\sum_{\gamma\in \Gamma} h^{p-F(\gamma);q-F(\gamma)}(H^{p+q-2F(\gamma)}(\M^\gamma;L_{\alpha})^\Gamma),$$ where $F(\gamma)$ is the fermionic shift, defined from the action of $\gamma$ on the tangent space of $\M$ at a $\gamma$-fixed point.

With these we can formulate our topological mirror symmetry conjecture:

\begin{conjecture}[\cite{hausel-thaddeus}] \label{tms} We have an agreement of Hodge numbers $h^{p,q}(\M)=h_{st,\alpha}^{p,q}(\M^\vee).$
\end{conjecture}
Note that as it stands the conjecture is about the Hodge numbers of the Higgs moduli spaces $\M=\M_{\Dol}$. However it was proved in  \cite{hausel-thaddeus} that the Hodge numbers of $\M_\Dol$ and $\M_{\DR}$ agree, and so \eqref{tms} is also about the agreement of (mixed) Hodge numbers of the proposed mirrors $\M_{\DR}$ and $\M^\vee_{\DR}$. This way we can interpret Conjecture~\ref{tms} as topological mirror symmetry for our SYZ mirror pair $\M_{\DR}$ and $\M^\vee_{\DR}$. 

Conjecture~\ref{tms} was proved for $\SL_2$ and $\SL_3$ in \cite{hausel-thaddeus}. More recently in 2020 the general case of $\SL_n$ was settled by Gr\"ochenig-Wyss-Ziegler \cite{groechenig-etal1}.  They used an arithmetic $p$-adic integration technique, pioneered  by Denef--Loeser \cite{denef-loeser} and used by Batyrev \cite{batyrev2} to check some topological mirror symmetry conjectures in the usual mirror symmetry. 

In \cite{hausel-gths} we observed a curious similarity. Namely Ng\^{o} in his proof \cite{ngo} of the Fundamental Lemma in the  Langlands program, reduced the Fundamental Lemma to the agreement of the number of points of certain singular Hitchin fibers over finite fields and in turn to the agreement of certain Hodge numbers of singular fibers of the Hitchin fibration. In \cite{hausel-gths}[\S 5.4] we argued that Ng\^{o}'s cohomological result is a relative version of the topological mirror symmetry conjecture  along the Hitchin fibration. A strategy was also proposed to deduce topological mirror symmetry using Ng\^{o}'s techniques. This proposal has been recently completed by Maulik--Shen in \cite{maulik-shen} in 2020 giving a new proof of the topological mirror symmetry Conjecture~\ref{tms} using Ng\^{o}'s techniques. 

Finally, in 2020,  Gr\"ochenig--Wyss--Ziegler in \cite{groechenig-etal2}  managed to extend their $p$-adic integration techniques from \cite{groechenig-etal1} for Higgs moduli spaces of general reductive groups $\G$ and in turn they found a new proof of Ng\^{o}'s cohomological result. 

\subsection{Geometric Langlands as enhanced Homological Mirror Symmetry} \label{gls}

In 2007 Kapustin--Witten \cite{kapustin-witten}  put forward a detailed circle of ideas amounting to  a physics derivation  of the geometric Langlands Correspondence as an enhanced mirror symmetry. They argued that a well-studied $S$-duality (or electro-magnetic or Montonen-Olive duality \cite{montonen-olive}) in a certain four-dimensional $N=4$ supersymmetric Yang-Mills theory, when reduced to two dimensions, yields an enhanced mirror symmetry, which in turn recovers the geometric Langlands correspondence as formulated by \cite{beilinson-drinfeld}. 

In this two-dimensional reduction, Montonen-Olive duality becomes an equivalence of a type $B$ sigma model with target the moduli space $\M_\DR$  of flat $\G$-connections on a complex curve $C$  and a type $A$ sigma model with target $\M^\vee_\DR$, the moduli space of flat $\G^\vee$-connections on $C$. As a consequence the category of boundary conditions in the two theories should be equivalent \beq \label{hmsgl} \mathscr{S}:D^b(\M_\DR)\simeq Fuk(\M_\DR^\vee),\eeq which can be interpreted as Kontsevich's homological mirror symmetry conjecture applied to the mirror pair $\M_\DR$ and $\M_\DR^\vee$ . 

Kapustin--Witten \cite{kapustin-witten} explained that this equivalence of categories has more structure due to the hyperk\"ahler targets, and more symmetries due to their origin in $4$-dimensions than the usual homological mirror symmetry, which arises from an equivalence of two $2$-dimensional sigma models. They use these additional ideas to construct from a flat connection $\alpha$ in $\M_\DR$, considered by its skycraper sheaf $\calO_\alpha\in D^b(\M_\DR)$, an element of $Fuk(\M^\vee_\DR)$ which they interpret as a $D$-module on the moduli space of $\G^\vee$-bundles on our curve $C$. The Hecke eigensheaf property  then in turn is deduced from the extra symmetry stemming from the $4$-dimensional origin.

First, due to the hyperk\"ahler targets, Kapustin--Witten talk about more structured branes (aka boundary conditions) by proposing that a brane should be either type $A$ or type $B$ with respect to all the three complex structures $I,J$ and $K$. This way they consider type $(B,A,A)$, $(A,B,A)$,$(A,A,B)$ and type $(B,B,B)$ branes on hyperk\"ahler manifolds. For example a type $(B,A,A)$ brane could be an $I$-holomorphic $\omega_\C=\omega_J+i\omega_K$ Lagrangian subvariety together with a local system. Or a type $(B,B,B)$ brane should be a hyperk\"ahler submanifold together with a hyperholomorphic connection on a bundle over it. 

In their framework \cite{kapustin-witten} argue that the mirror (S-dual) of an $(B,A,A)$ brane on the hyperk\"ahler $\M$ should be a $(B,B,B)$ brane on $\M^\vee$. In particular, if we just concentrate on complex structure $J$, that of $\M_{\DR}$, we see that the mirror of an $A$ brane should be a $B$ brane. 
The mirror relationships are slightly more subtle \cite[Table 2 p74]{kapustin-witten} in that the $A$-model in complex structure $J$ should be mirror to the $B$-model in complex structure $K$, which in turn, by \eqref{JK}, yields our version. 
This more refined version of mirror symmetry matches a type $B$ brane in complex structure $I$ to another type $B$ brane in complex structure $I$ on the mirror. 

This latter correspondence was also formulated by Donagi--Pantev \cite{donagi-pantev} as a classical limit - a first approximation of \eqref{hmsgl}- as \beq \label{classlim}\mathscr{S}:D^b(\M_\Dol)\simeq D^b(\M_\Dol^\vee)\eeq of the homological mirror symmetry \eqref{hmsgl}. First \cite{donagi-pantev} checks that \eqref{syzdol}  generically gives dual abelian varieties as fibers, for every reductive group $\G$. Second, they check that generically the Fourier-Mukai transform \cite{mukai} relative to the Hitchin base gives an equivalence like \eqref{classlim} which satisfies the additional intertwining of Hecke--Wilson symmetries discussed below. It is expected that \eqref{classlim} will have to be modified when extended over certain singular points of the moduli spaces.


Another direction of research - motivated by \cite{kapustin-witten}'s consideration of hyperk\"ahler branes - lead to new understandings of Lagrangian subvarieties  in $\M_\Dol$ and $\M_\DR$ and hyperholomorphic sheaves on $\M$. For example \cite{baraglia-schaposnik,garcia-prada-etal,hitchin-characteristic} studied various constructions of such hyperk\"ahler branes in all four different types, and contemplated  what their mirror should be. In particular, Hitchin \cite{hitchin-characteristic} proposed pairs of a $(B,A,A)$ brane on $\M$ and $(B,B,B)$ brane on $\M^\vee$ for which he could show that generically over the Hitchin base they are Fourier-Mukai dual. The $\G=\GL_2$ case of Hitchin's suggestion was the starting point of \cite{hausel-mellit-pei}, where an additional structure, the $\T$-action, played an important role. This point of view will be explained below in more detail. 

Second, due to the $4$-dimensional origin of their derivation of \eqref{hmsgl} Kapustin--Witten considered extra symmetries on these categories, arising from line operators in the $4$-dimensional theory. Namely, Wilson operators $\mathscr{W}^\mu$ attached to representations $\mu$ of $\G$ act on $D^b(\M_\DR)$ via tensoring with the vector bundle in the representation $\mu$ of the $\G$-bundle underlying the universal $\G$ flat connection.  On the other side Hecke operators (or t'Hooft operators for the physicists) $\mathscr{H}^\mu$ attached to irreducible representations $\mu$ of $\G$ act on the moduli space of $\G^\vee$ bundles and in turn on $D$-modules on them. The extra symmetry observation of \cite{kapustin-witten} is that these operators should intertwine the mirror symmetry of \eqref{hmsgl}. We will spell out these operators in  a more detailed way in the more symmetric classical limit \eqref{classlim} of \cite{donagi-pantev} in \S\ref{geosat}   below.  


The homological mirror symmetry \eqref{hmsgl} with these two additional structures: matching of hyperk\"ahler branes under mirror symmetry, and the Wilson-Hecke symmetry is what we call {\em enhanced mirror symmetry}. These go beyond the usual homological mirror symmetry of Kontsevich and stem from the supersymmetric and $4$-dimensional origins of $S$-duality.

\section{Enhanced mirror symmetry at the tip of the nilpotent cone}\label{enhanced}

The original motivation for the considerations below is to find a way to test the conjectured mirror pairs of $(B,A,A)$ and $(B,B,B)$ branes put forward in \cite{baraglia-schaposnik,garcia-prada-etal,hitchin-characteristic}. The only tests so far - which were often carried out in {\em loc.cit.} -  are to check if the proposed mirror pairs are indeed  Fourier-Mukai dual relative to the Hitchin maps. This can only be checked generically over the Hitchin base. We would like to see more global checks, in particular ones which can verify mirror symmetry proposals over the $0$-fibers of the Hitchin maps: the global nilpotent cones.  

We introduced a technique in \cite{hausel-mellit-pei,hausel-hitchin} which can verify mirror symmetry proposals over the nilpotent cone by considering the effect of mirror symmetry on morphisms in the corresponding categories. The difficulty to consider the morphisms in our categories arises from the non-compactness of our moduli spaces. For example the vector spaces of morphisms in the derived category $D(\M_\Dol)$ are typically infinite dimensional.  To measure their size, we will be looking at the $\T$-equivariant structure on them. Recall that the multiplicative group $\T:=\C^\times$ of the complex numbers acts on $\M$ by $\lambda:(E,\Phi)\to (E,\lambda\Phi)$, scalar multiplication of the Higgs field. Here we will be interested in a $\T$-equivariant extension of the classical limit \eqref{classlim} of the geometric Langlands correspondence, which should be as a first approximation an equivalence
$ \label{eqclasslim}\mathscr{S}:D_\T(\M_\Dol) \sim D_\T(\M^\vee_{\Dol})$ between the $\T$-equivariant derived categories of $\M_\Dol$ and $\M^\vee_{\Dol}$.  The morphisms between two objects $\calF_1$ and $\calF_2$ in $D_\T(\M_\Dol)$ can be identified with the graded vector space  $\Ext^*(X;\calF_1,\calF_2)$. To measure this graded vector space, we note that $\T$ acts on it, and assuming that the weight spaces are finite dimensional and vanish for large enough weights (which we expect for semiprojective varieties) 
we can define
the {\em equivariant Euler form} as \begin{multline*}\chi_\T(X;\calF_1,\calF_2)=\sum_{k,l} \dim(H^k({ R} {\mathcal Hom}(\calF_1,\calF_2))^l) (-1)^k t^{-l}=\\ =\sum_{k,l} \dim\left(\Hom_{{\mathbf D}_{coh}(X)}(\calF_1,\calF_2[k])^l\right) (-1)^k t^{-l}=\\=\sum_{k,l} \dim(\Ext^k(X;\calF_1,\calF_2)^l) (-1)^k t^{-l} \in \C((t)).\end{multline*}
With this we expect that $\mathscr{S}$ is an isometry: $$\chi_\T(\mathscr{S}(\calF),\mathscr{S}(\calG))=\chi_\T(\calF,\calG)$$

In \cite{hausel-mellit-pei} we  managed to check this isometry for several pairs of conjectured mirror branes from \cite{hitchin-characteristic}, while in \cite{hausel-hitchin} we checked this isometry for the conjectured mirror pairs relevant here; see below. In the second part of this paper we will recall some results of \cite{hausel-hitchin} about the mirror of very stable upward flows, introduce and study the multiplicity algebras of the Lagrangian upward flows following \cite{hausel-hitchin1}, and finally we will consider what the multiplicity algebra should correspond to on the mirror. 


\subsection{Very stable Higgs bundles and mirror symmetry}

The starting point is the recent paper \cite{hausel-hitchin}. First we recall its formalism. 

\subsubsection{Białynicki-Birula decomposition of semiprojective varieties}

Let the multiplicative group $\T:=\C^\times$ of the complex numbers act on a (possibly reducible) variety (a reduced separated  scheme of finite type over $\C$). We say that the action is {\em semiprojective} \cite[\S 1]{hausel-large} if the following three conditions hold \begin{enumerate} \item the action is {\em linear}, i.e. there is a locally closed $\T$-equivariant embedding of $X$ into $\P^N$ with a linear action of $\T$  (for example when $X$ is normal and quasi-projective)
	\item the fixed point subvariety $X^\T$ is proper and thus projective
	\item $\lim_{\lambda\to 0} \lambda\cdot x $ exists for all $x\in X$
\end{enumerate}
For $\alpha\in X^\T$ we define $W^+_\alpha:=\{x\in X| \lim_{\lambda\to 0} \lambda\cdot x = \alpha\}$ the {\em upward flow} from $\alpha$ and $W^-_\alpha:=\{x\in X| \lim_{\lambda\to \infty} \lambda\cdot x = \alpha\}$ the {\em downward flow} from $\alpha$. Then we have  $$X=\coprod_{\alpha\in X^\T} W_\alpha^+$$ the {Białyinicki-Birula partitio}n of $X$ and we define the projective variety $$\calC:=\coprod_{\alpha\in X^\T} W_\alpha^-$$ to be the {\em core }of $X$. We then have the following

\begin{theorem}$\!\!($\cite{bialynicki-birula},\cite[Proposition 2.1,Proposition 2.10]{hausel-hitchin}$)$. When $\alpha\in (X^s)^\T$, a $\T$-fixed point on the smooth locus, then $W^{\pm}_\alpha\subset X$ are locally closed subvarieties and $W^{\pm}_\alpha\cong T^{\pm}_\alpha X$ as $\T$-varieties. Moreover when $\omega\in \Omega^2(X^s)$ is a homogeneity $1$ symplectic form then $W^+_\alpha\subset X$ and $\calC\subset X$ are Lagrangian subvarieties. 
\end{theorem}

We call $\alpha\in (X^s)^\T$ and $W^+_\alpha$ {\em very stable} when $W^+_\alpha\cap \calC =\{\alpha\}$. Equivalently, $\alpha\in (X^s)^\T$ is very stable if and only if $W^+_\alpha\cap W_\beta^-\neq \emptyset $ implies $\beta=\alpha$. More generally one can show that the relation $\alpha\leq \beta$ when $W^+_\alpha\cap W_\beta^-\neq \emptyset$ induces a partial ordering. Then $\alpha\in (X^s)^\T$  is very stable if it is maximal with respect to this ordering. We then have

\begin{theorem}$\!\!($\cite[Proposition 2.14]{hausel-hitchin}.$)$\label{closed} $\alpha\in (X^s)^\T$ is very stable if and only $W_\alpha^+\subset X$ is closed. 
\end{theorem}

\subsubsection{Białynicki-Birula partition for  Higgs bundles} 

We will work with $\G=\GL_n$. We will denote by $\M$ the moduli space of semistable rank $n$ degree $d$ Higgs bundles $(E,\Phi)$ on a smooth projective curve of genus $g$. Here $E$ is a rank $n$ vector bundle of degree $d$ on $C$ and the Higgs field $\Phi\in H^0(C;\End(E)\otimes K)$. 

We have the {\em Hitchin map} \cite{hitchin-stable} given by the characteristic polynomial of the Higgs field: $$\begin{array}{cccc} h: &\M&\to& \calA:=\times_{i=1}^n H^0(C;K^i)\\ &(E,\Phi)&\mapsto& \det(x-\Phi).\end{array}$$ 

Then $h$ is a proper map \cite{hitchin,nitsure,simpson} and a completely integrable Hamiltonian system \cite{hitchin-stable,faltings} with respect to a natural holomorphic symplectic form on $\M$. In particular, the generic fibers are Lagrangian Abelian varieties, Jacobians of certain spectral curves. 

The $\T$-action on $\M$ is given by $(E,\Phi)\mapsto (E,\lambda\Phi)$. That makes the Hitchin map $\T$-equivariant if we let $\T$ act on $H^0(C;K^i)$ with weight $i$. As these weights are all positive on $\calA\cong \times_{i=1}^n H^0(C;K^i)$ it is semiprojective, and as the Hitchin map is proper and $\T$-equivariant we get that our $\T$-action on $\M$  is also  semiprojective. Additionally, we get that the core of $\calA$, i.e. the origin $0\in \calA$ pulls back to the core of $\M$. For $\M$ the core agrees with the nilpotent cone $h^{-1}(0)_{red}=\calC$. 

Generalising the notion of very stable bundle of Drinfeld and Laumon \cite{drinfeld2,laumon1} we can thus define a {\em very stable Higgs bundle} as a stable $\T$-fixed Higgs bundle $\calE\in {\M^s}^\T$ for which $W^+_\calE\cap h^{-1}(0)=\{\calE\}$ the only nilpotent Higgs bundle in its upward flow is itself. Thus by Theorem~\ref{closed} we know $\calE\in  {\M^s}^\T$ is very stable exactly when its upward flow is closed. To reformulate in terms of the Hitchin map we have an alternate version of Theorem~\ref{closed}.
\begin{theorem}[\cite{hausel-hitchin1} c.f. \cite{zelaci}] \label{proper}
	$\calE \in {\M^s}^\T$ is very stable if and only if $h^{-1}({0})\cap W^+_\calE=\{\calE\}$ if and only $h_\calE:=W_\calE^+\to \calA$ is proper if and only if it is finite. 
\end{theorem}

\begin{definition}\label{equmult} For $\calE\in \M^{s
		\T}$ define the rational function $$m_{\calE}(t):={{ \chi_\T(\Sym(T^{+*}_\calE))\over \chi_\T(\Sym(\calA^*))}}\in \Z(t).$$ 
	We call it the {\em equivariant multiplicity} of $W^+_\calE$. 
\end{definition}
We have the following
\begin{theorem}(\cite[Corollary 5.4]{hausel-hitchin}).\label{eqmul}
	When $\calE\in \M^{s
		\T}$ is very stable $m_{\calE}(t)$ is a polynomial \begin{itemize}  \item with non-negative coefficients, which is \item palindromic and \item monic, such that \item $m_{\calE}(1)=m_{F_\calE}$ is the multiplicity of the component $N_{F_\calE}\subset N$ in the nilpotent cone. \end{itemize}
\end{theorem}

Let $\ell\in \Z$, $m_i\in \N$ and $$\delta:=(\delta_0,\delta_1,\dots,\delta_{n-1})\in Jac_\ell(C)\times C^{[m_1]}\times \dots \times C^{[m_{n-1}]}$$ be a vector of representative divisors on $C$. To this we can construct a type $(1,\dots,1)$ $\T$-fixed Higgs bundle $\calE_\delta=(E_\delta,\Phi_\delta)$ where $$\calE_\delta=M_0\oplus\dots \oplus M_{n-1}$$ is a rank $n$ vector bundle $$M_i:=\calO(\delta_0+\dots+\delta_i)K^{-i}$$ and $$\Phi_\delta|_{M_i}:M_i\to M_{i+1}K\subset E_\delta K$$ is given by the defining section of $$H^0(C;M_i^{-1}M_{i+1}K)\cong H^0(C;\calO(\delta_i)).$$

The following classifies all very stable   $\T$-fixed type $(1,\dots,1)$ Higgs bundles and gives their equivariant multiplicity.

\begin{theorem}(\cite[Theorem 4.16,(5.18)]{hausel-hitchin}). \label{mult111} Let $\delta$ be as above and suppose that $\calE_\delta$ is a stable Higgs bundle. Then $\calE_\delta$ is very stable if and only if the effective divisor $\delta_1+\dots+\delta_{n-1}$ is reduced. Its equivariant multiplicity is given by $$m_{\calE_\delta}(t)=\prod_{i=1}^{n-1}\left[\begin{array}{c} n \\ i \end{array}\right]^{m_i}_t,$$  product of $t$-binomial coefficients.
\end{theorem}

\subsection{Multiplicity algebra and explicit Hitchin system on Lagrangians}

The main idea is to study, for $\calE\in \M^{s\T}$, the  restricted Hitchin map $$h_\calE:=h|_{W^+_\calE}: W^+_\calE\to \calA$$ in the framework of the Arnold school \cite{arnold-etal}. It is a $\T$-equivariant {\em Lagrangian map} between semiprojective vector spaces (i.e. only positive $\T$-weights) of the same dimension. Such maps are called {\em quasi-homogoneous} in \cite[\S 12.3]{arnold-etal}
. 

We recall from Theorem~\ref{proper} that $\calE$ is very stable if and only if  $h_\calE^{-1}(0)=\{\calE\}$ if and only if $h_\calE$ is proper.  For such maps - called {\em non-degenerate} in {\em loc.cit.}- \cite[\S 4,5,12]{arnold-etal} introduces and studies its local multiplicity algebra.

\begin{definition} When $\calE\in \M^{s\T}$ is very stable define $$Q_\calE:=Q_{h_\calE}:= \C[W^+_\calE]/(h_\calE^{-1}({\mathfrak m}_0))=\C[W^+_\calE]/(h_1,\dots,h_N)$$ the {\em local multiplicity algebra} of $h_\calE$ at $\calE$.  Here ${\mathfrak m}_0\subset \C[\calA]$ is the maximal ideal at $0\in \calA$ and  $$h_\calE=(h_1,\dots,h_N):\C^N\cong W^+_\calE\to \C^N\cong \calA$$ in some homogeneous coordinates.
\end{definition}

Scheme-theoretically $Q_\calE$ is just the coordinate ring of the scheme theoretical fiber of $h_\calE$ over $0$ or the scheme-theoretical intersection of $W^+_\calE\cap h^{-1}(0)$ of the upward flow with the nilpotent cone. Because $h_\calE$ is $\T$-equivariant we will get a $\T$-action, and thus a grading on $Q_\calE$. Because of this sometimes we call $Q_\calE$ the {\em equivariant multiplicity algebra} of $h_\calE$ at ${\calE}$. \begin{remark} Note that determining the algebra $Q_\calE$ explicitly by $N$ generators and $N$ relations gives us coordinates on $W^+_\calE$ such that the Hitchin map is given {\em explicitly} by the relations. \end{remark} Using results of \cite[\S 4,5,12]{arnold-etal}  we have the following

\begin{theorem}[\cite{hausel-hitchin1}] Let $\calE\in \M^s$ be very stable. Then its local multiplicity algebra is \begin{enumerate} \item finite dimensional,\item  graded $Q_\calE:=\bigoplus_{k=0}^m Q_\calE^k$ such that $Q_\calE^0\cong \C$, \item   Gorenstein, with socle $Q^m_\calE\cong \C J_{h_\calE}$, which is one-dimensional and spanned by the Jacobian $J_{h_\calE}$ of $h_\calE=(h_1,\dots,h_N)$ and  \item  a Poincar\'e duality ring. That is it has a natural bilinear pairing $(,):Q_\calE\times Q_\calE \to \C$ inducing a perfect pairing  $Q_\calE^k\times Q_\calE^{m-k}\to \C$ for all $k$. \item Finally, its Poincar\'e polynomial $\sum_{k=0}^m \dim(Q^k_\calE) t^k=m_\calE(t)$ agrees with the equivariant multiplicity. \end{enumerate}
\end{theorem}

This result gives a satisfactory explanation of all the properties of the equivariant multiplicity polynomial we observed in Theorem~\ref{eqmul}. The following then gives an explanation for the appearance of quantum binomial coefficients - which is well-known to be the Poincar\'e polynomial of a Grassmannian-  for the equivariant multiplicity. 

\begin{theorem}[\cite{hausel-hitchin1}] \label{mainema} Let $\delta=(\delta_0,\delta_1,\dots,\delta_{n-1})\in J_\ell(C)\times C^{[m_1]}\times \dots \times \times C^{[m_{n-1}]}$ be a vector of divisors. Assume that $\delta_1+\dots+\delta_{n-1}$ is reduced, then $\calE_\delta$ is very stable from Theorem~\ref{mult111} and its equivariant multiplicity algebra is $$Q_{\delta}\cong \prod_{i=1}^{n-1} H^*(\Gr(i,n);\C)^{m_i}$$ the product of the complex cohomology rings of Grassmanians $\Gr(i,n)$ of $i$-dimensional subspaces in $\C^n$.  
\end{theorem}

Our attack on Theorem~\ref{mainema} starts with understanding the map $h_{\delta}:W^+_{\delta}\to \calA$ using Hecke transformations. We explain it here in the first non-trivial case when $$\delta^k_c:=(\delta_0,\dots,\delta_{n-1})$$ where  $\delta_i=0$ unless $i=k$ and the divisor $\delta_k=c$ is one point.  For this we consider the $k$-th fundamental  Hecke correspondence over the Hitchin section $W^+_0:=W^+_{\calE_0}$ which coincides with the upward flow of the {\em canonical uniformising Higgs bundle} $\calE_0=(E_0,\Phi_0)$. Recall \cite[(3.16)]{hausel-hitchin} that it is given in chain notation as $$(E_0,\Phi_0)=(\calO\oplus K^{-1}\oplus \dots\oplus K^{1-n},\Phi_0)= \calO\stackrel{1}{\rightarrow} K^{-1} \stackrel{1}{\rightarrow} \dots \stackrel{1}{\rightarrow}  K^{1-n},$$ where morphisms on the arrows are twisted by $K$. The canonical uniformising Higgs bundle is at the top of the nilpotent cone. We will use Hecke transformations to generate all other type $(1,...,1)$ upward flows from the Hitchin section $W^+_0$.

We let \beq  \label{fundamental}\calH_k:=\left\{ (\calE_a,V)\in W^+_0\times \Gr(k, E_0|_{c})| \Phi_a|_{c}(V)\subset V\right\}\subset W^+_0\times \Gr(k, E_0|_{c}).\eeq It can be constructed \'etale locally over $\calA$ as the fixed point scheme of the self-map of $\Gr(k,E_0|_{c})$ induced by the invertible $\Phi_a|_{c}-\lambda I :E_0\to E_0$ where $\lambda\in \C$ is not an eigenvalue of $\Phi_a$. By performing a Hecke transformation of $\calE_a$ at $V$ for a point $(\calE_a,V)\in \calH_k$  we get that  $\calH_k\cong W^+_{\delta^k_c}$ and moreover  have the commutative diagram 

$$  \begin{array}{ccc}\calH_k&\stackrel{\cong}{\longrightarrow}& W^+_{\delta^k_c}\\ \!\!\!\!\!\!\!\!\! \pi_k{\Big\downarrow} && \,\,\,\,\,\,\, {\Big\downarrow} h_{{\delta^k_c}} \\ W^+_0 &\stackrel{h_0}{\longrightarrow}& \calA \end{array}$$

Hence $Q_{\delta^k_c}\cong Q_{\pi_k}$. In turn, the computation of $Q_{\pi_k}$  can be done in the Grassmannian $\Gr(k, E_0|_{c})$ and will yield in two different ways the two isomorphisms  \begin{multline} \label{cohgras} H^*(\Gr(k, E_0|_{c});\C)\cong Q_{\pi_k}\cong  \\\cong \frac{\C[p_1,\dots,p_k,q_1,\dots,q_{n-k}]}{\left((p_k+\dots+p_1x^{k-1}+x^k)(q_{n-k}+\dots+q_1x^{n-k-1}+x^{n-k})=x^n\right)}.\end{multline} Of course the isomorphism of the first and third rings gives the well-known presentation of the cohomology ring of the Grassmannian. 

In effect, we can think of the determination of the multiplicity algebra in \eqref{cohgras} to give coordinates on the upward flow $W^+_{\delta^k_c}$ so that the Hitchin map $h_{{\delta^k_c}}$ becomes explicit - basically given by the relations in the cohomology ring of the Grassmannian.

\subsection{Explicit Hitchin system for wobbly Lagrangians}

Here we show how one can generalise the technique above to understand multiplicity algebras of wobbly - i.e. not very stable- upward flows using the affine Grassmannian \cite{gortz,zhu}. We  start with a generalised notion of the $k$-th fundamental Hecke correspondence over the Hitchin section. Recall the {\em affine Grassmannian} $$\Gr:=\GL_n((z))/\GL_n[[z]],$$ where $\GL_n((z)):=\GL_n(\C((z)))$ and $\GL_n[[z]]:=\GL_n(\C[[z]])$.  It is a projective ind scheme, in particular its reduced is a nested union of projective varieties. It parametrises higher Hecke transformations of a vector bundle at a point on a curve.

Let $$\mu=(\mu_1\geq \mu_2\geq\dots \geq \mu_n)\in P^+\subset P\cong\Z^n$$ be a dominant weight and $$z^\mu:=\left(\begin{array}{cccc}z^{\mu_1} &0&\dots&0 \\ 0&z^{\mu_2}&\dots & 0\\ &&\vdots & \\ 0&0&\dots & z^{\mu_n}\end{array}\right)\in \GL_n(\C((z))).$$ We note that $\GL_n[[z]]$ acts from the left on $\Gr$ with orbit decomposition $$\Gr=\coprod_{\mu\in P^+} \Gr^\mu=\coprod_{\mu\in P^+} \GL_n[[z]] [z^\mu] $$ where $\Gr^\mu$ are labelled by dominant weights $\mu\in P^+$ as they are the orbits of $[z^\mu]$. We have a natural map \beq \label{fibr} \Gr^\mu\to \GL_n [z^\mu] \cong \GL_n/P_\mu \eeq given by setting $z=0$. We note that $\GL_n [z^\mu]\cong \GL_n/P_{\mu}$ is a partial flag variety and the map \eqref{fibr} is a finite rank vector bundle on $ \GL_n/P_{\mu}$. 
We denote by $\oGr^\mu$ the reduced of the closure of $\Gr^\mu$ in $\Gr$. We then have \beq\label{partition}\oGr^\mu=\coprod_{\mu \geq \lambda\in P^+}\Gr^\lambda\eeq where $\mu \geq \lambda$ is meant in the {\em dominance order}  on $P^+$i.e. when $\mu-\lambda$ is some sum of positive roots  (possibly with multiplicity). 

One important example is when $\mu=\omega_k=(\underbrace{1,\dots,1}_k,0,\dots,0)$ is the $k$-th fundamental weight. Then $\omega_k$ is minuscule (minimal in dominance order) and  $$\oGr^{\omega_k}=\Gr^{\omega_k}\cong \Gr(k,n)$$ the classical Grassmannian. 

For a dominant $\mu$ we can now define the Hecke correspondence of type $\mu$ over the Hitchin section as
\beq  \label{general}\calH^\mu:=\left\{ (\calE_a,V)\in W^+_0\times \Gr^\mu  \mbox{ s.t. }  \Phi_a|_{\Delta_c}(V)\subset V\right\}\subset W^+_0\times  \Gr^\mu,\eeq where $\Delta_c\to C$ is the formal neighbourhood of $c$ in $C$. In particular $\Delta_c\cong Spec(\C[[z]])$. We fix a trivialisation $E_0|_{\Delta_c}\cong \calO^n_{\Delta_c}$ and a trivialisation $K|_{\Delta_c}\cong \calO_{\Delta_c}$ then we can think of $\gamma_a:=\Phi_a|_{\Delta_c}\in \gl_n[[z]]$. Then if $V=[g]\in \Gr^\mu$ is represented by $g\in \GL_n((z))$  the condition $\Phi_a|_{\Delta_c}\in \gl_n[[z]]$ translates as \beq \label{affspring} \gamma_a g^{-1}\in \gl_n[[z]].\eeq For a fixed $a\in \calA$ this defines the {\em affine Springer fiber} of $\gamma_a$. 

By performing a type $\mu$ Hecke transformation of $\calE_0$ at $V\in \Gr^\mu$, this  will yield a new vector bundle together with a Higgs field, thanks to the invariance condition \eqref{affspring}. Provided some stability conditions are satisfied the new Higgs bundle will be stable and on the upward flow of $\calE_{\delta^\mu_c}$, where \beq \label{deltamuc}\delta^\mu_c=(\alpha_n c,\alpha_{n-1} c,\dots,\alpha_{1} c)\eeq by writing $\mu=\sum_{i} \alpha_i \omega_i$ in terms of the fundamental weights and $\alpha_i\in \N$ for $1\leq i \leq n-1$ and $\alpha_n\in \Z$. We claim that such a Hecke transformation will induce an isomorphism and so that we get the following diagram: \beq\label{diagram} \begin{array}{ccc}\calH^\mu&\stackrel{\cong}{\longrightarrow}& W^+_{\delta^\mu_c}\\ \!\!\!\!\!\!\!\!\! \pi_\mu{\Big\downarrow} && \,\,\,\,\,\,\, {\Big\downarrow} h_{{\delta^\mu_c}} \\ W^+_0 &\stackrel{h_0}{\longrightarrow}& \calA \end{array}\eeq
This way we  get $Q_{\pi_\mu}\cong Q_{h_{\delta^\mu_c}}$. Thus we  reduced the computation of the equivariant multiplicity algebra for a computation inside $\Gr^\mu$ by studying the equations describing $\calH^\mu\subset W^+_0 \times \Gr^\mu$. It turns out that this reduces to a relatively simple linear algebra computation. We have the following results and conjectures. In order to formulate them we will need to introduce the notion of dominant upward flows and their multiplicity. 

\begin{definition} Let $\calE\in \M^{s\T}$. We call the upward flow $W^+_\calE$ {\em dominant} if the Hitchin map $h_\calE:W^+_\calE\to \calA$ is dominant. In this case the induced map on algebra of functions $h_\calE^*: \C[{\calA}]\to \C[{W^+_\calE}]$ is injective and thus we get a extension $\C(\calA)\subset \C(W^+_\calE)$ of function fields. We define its degree to be the {\em multiplicity} of $W^+_\calE$  $$m_\calE:=[ \C(W^+_\calE):\C(\calA)].$$
\end{definition}

We note that $m_\calE$ also agrees with the rank $$m_\calE=\dim_{\C(\calA)} (M_\calE\otimes_{\C[\calA]} \C(\calA))$$ of the $\C[\calA]$-module $M_\calE$ given by $h_\calE^*$ and also with the cardinality of the generic fiber $$m_\calE=\#(h_\calE^{-1}(a)\cap W^+_\calE)$$ for generic $a\in\calA$. 
Notice that these last two definitions make sense for all upward flows. Dominance in turn then will be equivalent with nonzero multiplicity. 

We have the following results:

\begin{theorem}[\cite{hausel-hitchin1}] \label{conjecturewobbly} Let $\mu=(d+1)\omega_k\in P^+$  for $d\in \Z_{>0}$ and $1\leq k\leq n-1$ and $c\in C$.  When  $\calE_{\delta^\mu_c}\in \M^{s\T}$ the following hold
	\begin{enumerate}
		\item  Hecke modification of type $\mu$ of $W^+_0$  induces $\calH^\mu \to W^+_{\delta^\mu_c}$  an isomorphism
		\item $W^+_{\calE_{\delta^\mu_c}}$ is dominant
		\item $m_\calE=|W\cdot \mu|$ is the order of the Weyl orbit of $\mu$
		\item  we have
		$$Q_{\mu}\cong \C[J_d(Q_{\omega_k})]\cong \C\left[J_d(\Spec(H^*(\Gr(k,n),\C)))\right],$$ where for a scheme $X$ we denote by $J_d(X)$ the  $d-1$th jet scheme of $\spec(R)$.	In particular $J_d(X)(\C)=Hom(\Spec(\C[z]/(z^d)),X)$ is the set of $d-1$ jets in $X$. \end{enumerate}
\end{theorem}

\begin{remark} We conjecture that (1),(2) and (3) hold for any $\mu \in P^+$. 
\end{remark}

\begin{remark}
	
	As examples in the $n=2$ case let us give the multiplicity algebra for $d=1,2$ and $3$.
	First we have \beq \label{d=1} Q_{\delta_c^{2\omega_1}}\cong\C[a_0]/(a_0^2)\cong H^*(\P^1,\C).\eeq Then we have \beq \label{d=2} Q_{\delta_c^{2\omega_1}}\cong \C[a_0,a_1]/(a_0^2,a_0a_1).\eeq Note $(a_0^2,a_0a_1)=(a_0,a_1)^2\cap (a_0)$ thus the scheme theoretical intersection $Spec(Q_{\delta_c^{2\omega_1}})$ of the upward flow $W^+_{\delta_c^{2\omega_1}}$ and the nilpotent cone $h^{-1}(0)$ is the line $(a_0)$ with a double embedded point at the origin. Note that this upward flow was studied in \cite[\S 8.2]{hausel-hitchin}.  
	
	For $d=3$ we have the multiplicity $2$ algebra \beq \label{d=3} Q_{\delta_c^{3\omega_1}}\cong \C[a_0,a_1,a_2]/(a_0^2,a_0a_1,a_0a_2+a_1^2).\eeq Both \eqref{d=2} and \eqref{d=3} follows from Conjecture~\ref{conjecturewobbly}.4 and both can be proved by direct computation in $\Gr^\mu$ as explained above. 
\end{remark}

\begin{remark}It is surprising how complex $J_d(\Spec(H^*(\Gr(k,n),\C)))$ can be. In particular, in the $k=1$ case (i.e. jet schemes of the cohomology ring of projective space) there is only a conjecture about its multiplicity in \cite[Conjecture III.21]{yuenthesis}.
\end{remark}

\begin{remark} Finally we remark that already for type $(2)$ we have new phenomena. As discussed in \cite[\S 5.4]{hitchin4} there are multiplicity algebras depending on continuous parameters, in particular they cannot be isomorphic with cohomology rings, because cohomology rings are integral. \end{remark}

\subsubsection{Lagrangian closure of $W^+_\delta$}

\begin{definition} Let $\calE\in \M^{s\T}$. The {\em Lagrangian closure}  $\overline{\overline{W^+_\calE}}$ of 
	$W^+_\calE$ is the smallest closed union of upward flows containing $W^+_\calE$. In other words the Lagrangian closure is the closure in the quotient space by the BB partition. 
\end{definition}

Using \eqref{diagram} and \eqref{partition} we can deduce the following
\begin{theorem}[\cite{hausel-hitchin1}]  Let $\mu\in P^+$ and $c\in C$. Recall $\delta^\mu_c$ from \eqref{deltamuc}. Assume  $\calE_{\delta^\mu_c}\in \M^{s\T}$. Then 	 $$\overline{\overline{W^+_\calE}}=\coprod_{\mu\geq \lambda\in P^+} W^+_{\delta^\lambda_c},$$ i.e. the upward flows corresponding to dominant weights $\lambda$ less than or equal to $\mu$ in dominance order. 
\end{theorem}


\subsection{Towards a classical limit of the geometric Satake correspondence}
\label{geosat}
Finally, we will formulate some conjectures which were the original motivation of much of the previous ideas. In particular, they hint at a new construction of the irreducible representations of $\GL_n(\C)$, and  more generally of any complex reductive group $\G$. 

The general setup comes from the classical limit \eqref{classlim} of the geometric Langlands program, as formulated 
in \cite{donagi-pantev}. Here we sketch some of the expectations of this classical limit in a schematic (not completely well-defined) manner.  It should be an equivalence of some sort of derived categories of coherent sheaves $$\mathscr{S}:D^b(\M_\Dol)\to D^b(\M_\Dol^\vee).$$ Several properties of this equivalence were proposed and some established in \cite{donagi-pantev}. In particular, $\mathscr{S}$ should be a relative Fourier-Mukai transform along the generic locus in $\calA_\G\cong \calA_{\G^\vee}$. Another  crucial property \cite{kapustin-witten} - which we called enhanced mirror symmetry in \S \ref{gls} above-  is that $\mathscr{S}$ should intertwine
the actions of certain Hecke operators on $D^b(\M_\Dol)$ and the Wilson operators on $D^b(\M_\Dol)$. Let $\mu\in X_+(\G^\vee)=X^+(\G)$ be a dominant character of $\G^\vee$. We denote by $$\calH^\mu:=\{(E,\Phi)\in \M_\Dol, [g]\in \Gr^\mu | g^{-1}\Phi_c g \in \G[[z]]\}\subset \M_\Dol\times \Gr^\mu$$ some space of Hecke correspondences at a point $c\in C$. Indeed this gives us $$\begin{array}{ccc} &\calH^\mu&\\ \pi_\mu\swarrow&&\searrow f^\mu \\\!\!\!\!\!\! \M_\Dol&&\,\,\,\,\,\, \M_\Dol \end{array}$$ two maps to $\M_\Dol$, first the projection $\pi_\mu$ to the first factor, and second $f^\mu$, the Hecke transformation\footnote{Here we ignore issues about stability.} of $(E,\Phi)$ by the compatible Hecke transform $[g]\in \Gr^\mu$,  which is expected to induce $$\mathscr{H}^\mu :=f^\mu_* \circ \pi_\mu^*: D^b(\M_\Dol)\to  D^b(\M_\Dol)$$ a Hecke or the physicists' {\em t'Hooft operator}. 

On the other hand we can consider the so-called {\em Wilson} operators $$\mathscr{W}^\mu : \begin{array}{ccc} D^b(\M^\vee_{Dol}) &\to& D^b(\M^\vee_\Dol) \\ \calF &\mapsto& \calF\otimes \rho_\mu(\bE|_{\M^\vee_{\Dol}\times \{c\}}) \end{array} $$ given by tensoring with the  universal $\G^\vee$ bundle $\bE$ in the representation $\rho_\mu: \G^\vee\to \GL(V_{\rho_\mu})$.

We then expect \cite{donagi-pantev,kapustin-witten} that \beq\label{intertwine}\mathscr{W}^\mu \circ \mathscr{S} = \mathscr{S} \circ \mathscr{H}^\mu.\eeq

There are two more expectations for the classical limit $\mathscr{S}$, both are motivated from Fourier-Mukai transform where the analogous statements hold. First we  expect that for any $\calF\in D^b(\M_\Dol)$ we should have \beq \label{Spush} (h_\G)_*(\calF)\cong \mathscr{S}(\calF)|_{W^+_0}.\eeq Second, the structure sheaf of the Hitchin section should be mirror to the structure sheaf of the mirror Higgs moduli space: \beq \label{Sstructure}\mathscr{S}(\calO_{W^+_0})\cong \calO_{\M^\vee_{\Dol}} .\eeq

Combining \eqref{intertwine} with \eqref{Sstructure} we can deduce that $$\mathscr{S}(\mathscr{H}^\mu(\calO_{W^+_0}))=\mathscr{W}^\mu(\mathscr{S}(\calO_{W^+_0}))=\mathscr{W}^\mu(\calO_{\M^\vee_{\Dol}}).$$  On the one hand we should have $${\mathrm supp}(\mathscr{H}^\mu(\calO_{W^+_0}))=\overline{\overline{W^+_\mu}},$$ where $W_\mu^+$ is the upward flow from a certain $\calE_\mu$ maximally split $G$-Higgs bundle of type $\mu$ at $c\in C$. On the other hand $$\mathscr{W}^\mu(\calO_{\M^\vee_{\Dol}})=\rho_\mu(\bE)_c=:\Lambda_\mu,$$ the vector bundle associated  to the principal bundle $\bE_c$ in the representation $\rho_\mu$. 

Thus Kapustin-Witten's \eqref{intertwine} implies $$\mathscr{S}(\mathscr{H}^\mu(\calO_{W^+_0}))=\Lambda_\mu.$$ We can test this by \eqref{Spush}  $$\Lambda_\mu|_{ ^{L}W_0^+}= \mathscr{S}(\mathscr{H}^\mu(\calO_{W^+_0}))|_{ ^{L}W_0^+}=(h_\G)_*(\mathscr{H}^\mu(\calO_{W^+_0})).$$

In \cite{hausel-hitchin} we have made arguments that the mirror of the structure sheaf of a very stable type $1,\dots,1$ upward flow $W_{\delta}$ is $$\Lambda_\delta:=\bigotimes_{i=1}^{n-1} \bigotimes_{j=1}^{m_i} \Lambda^i\E_{c_{ij}},$$ where $(\E,\bPhi)$ is a universal Higgs bundle on $\M\times C$ and $$\delta_i=c_{i1}+{c_{i2}}+\dots \in C^{[m_i]}.$$ In particular, one expectation of mirror symmetry
is that $$h_*(\calO_{W^+_\delta})\cong \Lambda_\delta|_{W^+_0}.$$ This follows from Theorem~\ref{mult111} and a direct computation for $\chi_\T (\E_c)$.

In \cite[\S 8.2]{hausel-hitchin} we proposed that for $n=2$ the mirror of $Sym^2(\E_c)$ should be the structure sheaf of the Lagrangian closure $\overline{\overline{W^+_{\delta_{c}^2}}}$ where $\delta_{c}^2=(0,2c)$. We can generalise this as
follows. 

\begin{conjecture}\label{geomsati} Let $c\in C$ and $\G$ a reductive group. Then we have the following conjectures. \begin{enumerate}
		\item For any $\mu\in X_+(\G^\vee)$ the support of the mirror of  $\rho_\mu(\bE_c)$ is $\overline{\overline{W^+_{\delta_c^\mu}}}$
		\item Let $\mu\in X_+(\G^\vee)$ such that the corresponding irreducible $\G^\vee$ representation $
		\rho_\mu$ is multiplicity free. Then the mirror of $\rho_\mu(\bE_c)$ is $\calO_{\overline{\overline{W^+_{\delta_c^\mu}}}}$. 
		\item In the latter case the multiplicity algebra of the restricted Hitchin map $h_\G:\overline{\overline{W^+_{\delta_c^\mu}}}\to \calA$ is isomorphic with the cohomology ring of $\overline{\Gr}^\mu$. 
	\end{enumerate}
\end{conjecture}

\begin{remark} In \cite{feigin-frenkel-rybnikov}, studying opers in the geometric Langland program, the authors construct a canonical Poincar\'e duality ring structure on the underlying vector space $V_\mu$ of all irreducible representation $\rho_\mu$ of $\G^\vee$. In the case when $\rho_\mu$ is multiplicity-free this ring is isomorphic with the cohomology ring $H^*(\overline{\Gr}^\mu)$. Note that according to 
	\cite[Theorem 1.5]{ginzburg} these are precisely the cases when $$H^*(\overline{\Gr}^\mu)\cong  I\!H^*(\overline{\Gr}^\mu),$$ when the cohomology ring satisfies Poincar\'e duality. In this case this ring was more carefully studied in \cite{panyushev}.
\end{remark}


{\noindent \bf Acknowledgements.} 
The author thanks Nigel Hitchin for introducting him to Higgs bundles during 1995-1998, suggesting the SYZ picture for Langlands dual Hitchin systems in 1996 and for the more recent collaborations \cite{hausel-hitchin,hausel-hitchin1}. He also thanks David Ben-Zvi, Pierre-Henri Chaudouard, Pierre Deligne, Ron Donagi, Sergei Gukov, Jochen Heinloth, Joel Kamnitzer, G\'erard Laumon, Anton Mellit, David Nadler, Andy Neitzke,  Ng\^{o} Bao Ch\^au, Michael Thaddeus, Tony Pantev, Du Pei, Rich\'ard Rim\'anyi, Leonid Rybnikov, Vivek Shende, Bal\'azs Szendrői, Andr\'as Szenes, Fernando Rodriguez-Villegas, Edward Witten and Zhiwei Yun for discussions about the subjects in this paper over the years. Thanks are also due to H\"ulya Arg\"uz, Jakub L\"owit, Bal\'azs Szendr{\H{o}}i and Nigel Hitchin for careful reading of the paper.


\end{document}